\documentclass[11pt,reqno]{amsart}
\usepackage[T1]{fontenc}
\usepackage[utf8]{inputenc}
\usepackage{lmodern}
\usepackage{amsmath,amssymb,amsthm,mathtools}
\usepackage{microtype}
\usepackage[margin=1.05in]{geometry}
\usepackage[hidelinks]{hyperref}
\newtheorem{theorem}{Theorem}[section]
\newtheorem{lemma}[theorem]{Lemma}

\newtheorem{conjecture}[theorem]{Conjecture}
\theoremstyle{remark}

\numberwithin{equation}{section}
\newcommand{\R}{\mathbb R}
\newcommand{\tr}{\operatorname{tr}}
\newcommand{\op}{\mathrm{op}}
\newcommand{\Frob}{\mathrm F}
\newcommand{\one}{\mathbf 1}
\allowdisplaybreaks[0]
\makeatletter
\let\original@setauthors\@setauthors
\renewcommand{\@setauthors}{%
  \begingroup
  \let\MakeUppercase\@firstofone
  \def\author@andify{%
    \nxandlist
      {\unskip,\penalty-1\space\ignorespaces}%
      {\unskip,\penalty-2\space\ignorespaces}%
      {\unskip,\penalty-2\space\ignorespaces}%
  }%
  \original@setauthors
  \endgroup
}
\makeatother

\title[Triangles without a common vertex]
{Spectral supersaturation for triangles without a common vertex}
\author[Dongxiu Cai]{Dongxiu Cai\textsuperscript{1}}
\author[Xiao-Dong Zhang]{Xiao-Dong Zhang\textsuperscript{1,2}}
\thanks{\textsuperscript{1} School of Mathematical Sciences, MOE-LSC,
SHL-MAC, Shanghai Jiao Tong University, 800 Dongchuan Road,
Shanghai 200240, P.R.~China. Email:
\texttt{diudiutse@sjtu.edu.cn} (Dongxiu Cai).}
\thanks{\textsuperscript{2} Corresponding author: Xiao-Dong Zhang.
Email: \texttt{xiaodong@sjtu.edu.cn}.}
\date{}
\hypersetup{pdftitle={Spectral supersaturation for triangles without a common vertex},
  pdfauthor={Dongxiu Cai and Xiao-Dong Zhang},
  pdfsubject={A proof of the Li--Feng--Peng triangle counting conjecture},
  pdfkeywords={spectral supersaturation, triangles, spectral radius}}

\begin{document}
\begin{abstract}
Let $G$ be a graph with $m$ edges, $t(G)$ triangles, and adjacency
spectral radius $\lambda(G)>\sqrt m$. Li, Feng and Peng (2025)
conjectured that $t(G)\ge\sqrt m-O(1)$ if no vertex belongs to
every triangle of $G$. We prove the stronger bound
$t(G)\ge\lambda(G)-C$ for an absolute constant $C$, and give an
explicit choice of $C$. The proof bounds the sum of the squared
adjacency eigenvalues other than the largest and smallest.
Interlacing and matrix perturbation estimates then restrict the
neighbourhood of a vertex corresponding to a largest entry of a
Perron vector to two disjoint edges and isolated vertices.
A weighted count of the triangles completes the proof.
\end{abstract}
\maketitle
\raggedbottom

\section{Introduction}
\label{sec:introduction}

Let $G$ be a finite simple graph with $n$ vertices and $m$ edges.
Write its adjacency eigenvalues in decreasing order as
\[
 \lambda_1(G)\ge\lambda_2(G)\ge\cdots\ge\lambda_n(G),
\]
and let $\lambda(G)=\lambda_1(G)$ be its adjacency spectral radius.
Write $t(G)$ for its number of triangles, and put
\[
 S(G)=\sum_{i=2}^{n-1}\lambda_i(G)^2,
\]
where an empty sum is understood to be zero.
Nosal's spectral version
of Mantel's theorem~\cite{Nosal1970} states that
$\lambda(G)>\sqrt m$ forces a triangle. The threshold is attained
by every complete bipartite graph. Quantitative refinements ask
how many triangles must occur when the spectral radius reaches
or exceeds this threshold. Bollob\'as and
Nikiforov~\cite{BollobasNikiforov2007} proved that
\begin{equation}
\label{eq:intro-moment-bound}
 t(G)\ge\frac{\lambda(G)\bigl(\lambda(G)^2-m\bigr)}3.
\end{equation}
Ning and Zhai~\cite{NingZhai2023} established the sharp bound
\begin{equation}
\label{eq:intro-NZ}
 t(G)\ge\left\lfloor\frac{\sqrt m-1}{2}\right\rfloor
\end{equation}
for graphs with $\lambda(G)\ge\sqrt m$, unless $G$ is a complete
bipartite graph, possibly together with isolated vertices.
Note that equality in \eqref{eq:intro-NZ} is attained by adding
one edge inside the larger part of $K_{4b+4,b}$, for an integer
$b\ge1$; all triangles contain the added edge.

For edge conditions, Xiao and
Katona~\cite{XiaoKatona2021} showed that an $n$-vertex graph
with more than $\lfloor n^2/4\rfloor$ edges has at least $n-2$
triangles if no vertex belongs to every triangle.
Li, Feng and Peng~\cite[Theorem~2.2]{LiFengPeng2025} proved a
spectral analogue: if $n\ge113$,
$\lambda(G)\ge\sqrt{\lfloor n^2/4\rfloor}$, and the triangle
covering number of $G$ is at least $2$, then $t(G)\ge n-3$.
The triangle covering number is the minimum size of a vertex set
that meets every triangle.
Motivated by \eqref{eq:intro-NZ}, they proposed the following
conjecture.

\begin{conjecture}[{Li--Feng--Peng~\cite[Conjecture~6.1]{LiFengPeng2025}}]
\label{conj:LFP}
If $G$ is an $m$-edge graph with $\lambda(G)>\sqrt m$, and there
is no vertex contained in all triangles of $G$, then $G$ contains
at least $\sqrt m-O(1)$ triangles.
\end{conjecture}

In this paper, we prove a stronger version of
Conjecture~\ref{conj:LFP}.

\begin{theorem}
\label{thm:main}
Put
\[
 C=\binom{100000005}{6}+10.
\]
Let $G$ be a simple graph on $n$ vertices with $m$ edges satisfying
$\lambda=\lambda(G)>\sqrt m$. If no vertex of $G$ belongs to every
triangle, then
\begin{equation}
\label{eq:main}
 t(G)\ge\lambda(G)-C>\sqrt m-C.
\end{equation}
In particular, Conjecture~\ref{conj:LFP} holds.
\end{theorem}

The coefficient of $\sqrt m$ in \eqref{eq:main} is best possible.
The construction
of Li, Feng and Peng~\cite{LiFengPeng2025} obtained by adding two
disjoint edges inside the larger part of $K_{4b+3,b}$ has
$m=4b^2+3b+2$ and $t(G)=2b$ for $b\ge2$. It satisfies the
hypotheses, and $t(G)/\sqrt m$ tends to $1$.

Several recent results give more triangles under a stronger
spectral assumption. Zhang and Zhai~\cite{ZhangZhai2026} proved,
for sufficiently large $m$, that
\[
 \lambda(G)\ge\frac{1+\sqrt{4m-3}}2
 \quad\Longrightarrow\quad t(G)\ge\frac{m-1}{2}.
\]
Chen, Li and Tang~\cite[Theorem~1.3]{ChenLiTang2026} showed that
every $m$-edge graph satisfies
\[
 t(G)\ge m\bigl(\lambda(G)-\sqrt m\bigr).
\]
Theorem~\ref{thm:main} gives a uniform lower bound when the
spectral surplus $\lambda(G)-\sqrt m$ is arbitrarily small,
using the condition that the triangles have no common vertex.

We briefly describe the proof. Suppose that $G$ violates
\eqref{eq:main}. The second and third spectral moments imply that
$S(G)<41/10$. Interlacing passes this bound to every induced
subgraph. We use a second-order perturbation estimate to show
that $S(J)>9/2$ for two types of induced subgraphs $J$.
Their exclusion, together with the classical Ramsey bound,
determines the neighbourhood of a vertex corresponding to a
largest entry of a nonnegative Perron vector. The neighbours span
exactly two disjoint edges, and the nonneighbours form an
independent set. Finally, the four eigenvector equations at the
endpoints of those two edges yield
$\lambda(G)<t(G)+2$, contradicting the assumed violation.

The article is organized as follows. Section~\ref{sec:conventions}
sets out the notation and spectral moment identities.
In Section~\ref{sec:middle}, we bound the middle eigenvalues and
prove the perturbation estimates used for induced subgraphs.
Section~\ref{sec:configurations} applies these estimates to the
two configurations needed in the structural argument.
In Section~\ref{sec:neighbourhood}, the Perron equations and a
Ramsey bound give the neighbourhood structure of a putative
counterexample. Section~\ref{sec:main-proof} proves the weighted
triangle estimate and completes the proof of
Theorem~\ref{thm:main}.

\section{Conventions}
\label{sec:conventions}

\subsection{Graphs and spectral notation}

All graphs are finite, simple, and undirected; they need not be
connected. For a graph $G$, let $A_G$ be its adjacency matrix and write
\[
 n=|V(G)|,\qquad m=e(G)=|E(G)|,\qquad t=t(G).
\]
Here $t(G)$ counts all three-element vertex sets spanning a triangle,
including triangles contained in larger cliques. When $G$ is fixed,
we write $A=A_G$. For $n\ge1$, arrange the adjacency eigenvalues,
with multiplicity, as
\[
 \lambda_1(G)\ge\lambda_2(G)\ge\cdots\ge\lambda_n(G).
\]
By the Perron--Frobenius theorem (see~\cite{BrouwerHaemers2012}),
the adjacency spectral radius is
$\lambda(G)=\lambda_1(G)$. When the graph is fixed, write
$\lambda=\lambda(G)$ and $\lambda_i=\lambda_i(G)$; thus
$-\lambda\le\lambda_i\le\lambda$ for every $i$. The empty graph has
spectral radius zero. Graphs with at most one vertex do not satisfy
$\lambda>\sqrt m$.

For $Z\subseteq V(G)$, let $G[Z]$ be the induced subgraph and write
\[
 N(v)=N_G(v),\qquad d_Z(v)=|N(v)\cap Z|.
\]
Adjacency and neighbourhoods are always taken in the graph currently
under discussion. A graph has \emph{no common triangle vertex} if no
vertex belongs to every triangle, or equivalently, if for each vertex
$v$ there is a triangle avoiding $v$.

The two spectral moment identities used below are
\begin{equation}
\label{eq:trace2}
 \sum_{i=1}^n\lambda_i^2=\tr(A^2)=2m,
\end{equation}
and
\begin{equation}
\label{eq:trace3}
 \sum_{i=1}^n\lambda_i^3=\tr(A^3)=6t.
\end{equation}
The first identity counts the two orientations of each edge; the
second counts the three starting vertices and two directions around
each triangle. Also $\tr A=0$, so $\lambda_n\le0$.

\subsection{Matrices, coordinates, and the middle eigenvalues}

For $Z\subseteq V(G)$, set
\[
 A_Z=A_{G[Z]}.
\]
Vectors are columns. Their norm is the Euclidean norm, and
$\|M\|_{\op}$ is the induced operator norm. The Frobenius norm is
\[
 \|M\|_{\Frob}^2=\sum_{i,j}M_{ij}^2.
\]
We use $\|M\|_{\op}\le\|M\|_{\Frob}$, which follows from
Cauchy--Schwarz; for symmetric $M$, $\|M\|_{\Frob}^2=\tr(M^2)$.
In block matrices and vectors, $\one$ denotes an all-ones column
vector of the size determined by its block.
Block matrices below always come with an
explicit order of their coordinate sets.

For a real symmetric matrix $M$ of order $r\ge2$, use the analogous
notation
\begin{equation}
\label{eq:middle-mass}
 S(M)=\sum_{i=2}^{r-1}\lambda_i(M)^2,
 \qquad S(A_G)=S(G).
\end{equation}
The sum is zero when $r=2$.

\section{Middle eigenvalues and matrix estimates}
\label{sec:middle}

We first bound the middle eigenvalues of a graph with few triangles.
The perturbation estimates at the end of this section will give
lower bounds on the same quantity for the induced subgraphs used
in Section~\ref{sec:configurations}.

\begin{lemma}[Middle eigenvalues of a principal submatrix]
\label{lem:interlacing}
Let $B$ be a principal submatrix of order $r\ge2$ of a real
symmetric matrix $M$ of order $n$. Then $S(B)\le S(M)$.
Consequently, a graph satisfying $S(G)<41/10$ contains no $K_7$.
\end{lemma}

\begin{proof}
By the Cauchy interlacing theorem,
\begin{equation}
\label{eq:interlacing}
 \lambda_j(M)\ge\lambda_j(B)\ge\lambda_{n-r+j}(M)
 \qquad(1\le j\le r).
\end{equation}
For $2\le j\le r-1$, use the upper bound if $\lambda_j(B)\ge0$
and the lower bound if $\lambda_j(B)<0$. This gives
\[
 \lambda_j(B)^2\le
 \begin{cases}
  \lambda_j(M)^2,&\lambda_j(B)\ge0,\\
  \lambda_{n-r+j}(M)^2,&\lambda_j(B)<0.
 \end{cases}
\]
All indices on the right lie in $\{2,\ldots,n-1\}$ and are
distinct within each case. They are also distinct between cases:
if $\lambda_j(B)\ge0>\lambda_k(B)$, then
$j<k\le n-r+k$. Summing proves $S(B)\le S(M)$.

The adjacency matrix of $K_7$ is $J_7-I_7$. Its eigenvalues are
$6$ and six copies of $-1$, so $S(K_7)=5$. A seven-vertex clique
would therefore give $5\le S(G)$ by the first conclusion, which
is impossible when $S(G)<41/10$. This finishes the proof.
\end{proof}

\begin{lemma}[A bound from the number of triangles]
\label{lem:smallmass}
Let $G$ be a simple graph on $n$ vertices with $m$ edges satisfying
$\lambda=\lambda(G)>\sqrt m$, $t=t(G)<\lambda$, and
$\lambda\ge1000$. Then $S(G)<41/10$.
\end{lemma}

\begin{proof}
If $n=2$, then $S(G)=0$, so assume $n\ge3$. Combining the two
spectral moments gives
\begin{equation}
\label{eq:nonnegative-moment}
\begin{alignedat}{2}
 6t&=\lambda^3+\sum_{i=2}^n\lambda_i^3
 &&\quad\text{(by \eqref{eq:trace3})}\\
 &=\lambda^3+\sum_{i=2}^n(\lambda+\lambda_i)\lambda_i^2
            -\lambda\sum_{i=2}^n\lambda_i^2\\
 &=\lambda^3+\sum_{i=2}^n(\lambda+\lambda_i)\lambda_i^2
            -\lambda(2m-\lambda^2)
 &&\quad\text{(by \eqref{eq:trace2})}\\
 &=\sum_{i=2}^n(\lambda+\lambda_i)\lambda_i^2
              +2\lambda(\lambda^2-m).
\end{alignedat}
\end{equation}
Since $\lambda_i\ge-\lambda$, every term on the right is
nonnegative. Hence
\[
\begin{alignedat}{2}
 0<\lambda^2-m&\le\frac{6t}{2\lambda}
 &&\quad\text{(by \eqref{eq:nonnegative-moment})}\\
 &=\frac{3t}{\lambda}\\
 &<3
 &&\quad\text{(by $t<\lambda$)}.
\end{alignedat}
\]

We next show that $\lambda+\lambda_i>\lambda/4$ for every middle
eigenvalue. If all middle eigenvalues are nonnegative, then for
every $2\le i\le n-1$,
\[
\begin{alignedat}{2}
 \lambda+\lambda_i&\ge\lambda
 &&\quad\text{(by $\lambda_i\ge0$)}\\
 &>\lambda/4
 &&\quad\text{(by $\lambda>0$)}.
\end{alignedat}
\]
If there is a negative middle eigenvalue, take any index
$2\le i\le n-1$ with $\lambda_i<0$. Then
\[
\begin{alignedat}{2}
 2\lambda_i^2&\le\lambda_i^2+\lambda_n^2
 &&\quad\text{(by $\lambda_n\le\lambda_i<0$)}\\
 &\le\sum_{j=2}^n\lambda_j^2\\
 &=2m-\lambda^2
 &&\quad\text{(by \eqref{eq:trace2})}\\
 &<\lambda^2
 &&\quad\text{(by $m<\lambda^2$)}.
\end{alignedat}
\]
The inequality $\lambda_i\le-3\lambda/4$ would instead give
$2\lambda_i^2\ge9\lambda^2/8>\lambda^2$. Therefore
\[
\begin{alignedat}{2}
 \lambda+\lambda_i&>\lambda-3\lambda/4
 &&\quad\text{(by $\lambda_i>-3\lambda/4$)}\\
 &=\lambda/4.
\end{alignedat}
\]
Every nonnegative middle eigenvalue in this second case satisfies
the estimate from the first case. Thus the bound holds for all
middle eigenvalues, whether or not a negative one exists.

It follows that
\begin{equation}
\label{eq:middle-mass-24}
\begin{alignedat}{2}
 S(G)&\le\frac4\lambda
          \sum_{i=2}^{n-1}(\lambda+\lambda_i)\lambda_i^2
 &&\quad\text{(by $\lambda+\lambda_i>\lambda/4$)}\\
 &\le\frac4\lambda\sum_{i=2}^n(\lambda+\lambda_i)\lambda_i^2\\
 &=\frac{24t}{\lambda}-8(\lambda^2-m)
 &&\quad\text{(by \eqref{eq:nonnegative-moment})}\\
 &<24
 &&\quad\text{(by $t<\lambda$, $m<\lambda^2$)}.
\end{alignedat}
\end{equation}
In particular, $\lambda_i^2\le S(G)<25$ and $\lambda_i>-5$ for
$2\le i\le n-1$. Hence
\begin{equation}
\label{eq:middle-cube}
\begin{alignedat}{2}
 \sum_{i=2}^{n-1}\lambda_i^3
 &=\sum_{i=2}^{n-1}\lambda_i\lambda_i^2\\
 &\ge-5\sum_{i=2}^{n-1}\lambda_i^2
 &&\quad\text{(by $\lambda_i>-5$)}\\
 &=-5S(G)
 &&\quad\text{(by \eqref{eq:middle-mass})}.
\end{alignedat}
\end{equation}
Also,
\begin{equation}
\label{eq:extreme-square-difference}
\begin{alignedat}{2}
 0<\lambda^2-\lambda_n^2&=2(\lambda^2-m)+S(G)
 &&\quad\text{(by \eqref{eq:trace2})}\\
 &<2\cdot3+24
 &&\quad\text{(by the bounds just proved)}\\
 &=30.
\end{alignedat}
\end{equation}

The following factorization relates the extreme cubes to the
middle squares:
\begin{equation}
\label{eq:extreme-factorization}
\begin{alignedat}{2}
 &2\lambda(\lambda^3+\lambda_n^3)
   -3\lambda^2(\lambda^2-\lambda_n^2)
   +3(\lambda^2-\lambda_n^2)^2\\
 &=(\lambda+\lambda_n)^2
   (2\lambda^2-4\lambda\lambda_n+3\lambda_n^2)\\
 &\ge0
 &&\quad\text{(by $\lambda_n\le0<\lambda$)}.
\end{alignedat}
\end{equation}
Consequently,
\begin{equation}
\label{eq:S-cubic-bound}
\begin{alignedat}{2}
 3\lambda^2S(G)
 &=3\lambda^2[\lambda^2-\lambda_n^2-2(\lambda^2-m)]
 &&\quad\text{(by \eqref{eq:trace2})}\\
 &=3\lambda^2(\lambda^2-\lambda_n^2)-6\lambda^2(\lambda^2-m)\\
 &\le3\lambda^2(\lambda^2-\lambda_n^2)
 &&\quad\text{(by $\lambda^2>m$)}\\
 &\le2\lambda(\lambda^3+\lambda_n^3)
       +3(\lambda^2-\lambda_n^2)^2
 &&\quad\text{(by \eqref{eq:extreme-factorization})}.
\end{alignedat}
\end{equation}
We estimate the right-hand side one term at a time:
\[
\begin{alignedat}{2}
 3\lambda^2S(G)
 &\le2\lambda(\lambda^3+\lambda_n^3)
          +3(\lambda^2-\lambda_n^2)^2
 &&\quad\text{(by \eqref{eq:S-cubic-bound})}\\
 &=2\lambda\left(6t-\sum_{i=2}^{n-1}\lambda_i^3\right)
          +3(\lambda^2-\lambda_n^2)^2
 &&\quad\text{(by \eqref{eq:trace3})}\\
 &\le2\lambda[6t+5S(G)]+3(\lambda^2-\lambda_n^2)^2
 &&\quad\text{(by \eqref{eq:middle-cube})}\\
 &<2\lambda[6\lambda+5S(G)]+3(\lambda^2-\lambda_n^2)^2
 &&\quad\text{(by $t<\lambda$)}\\
 &<2\lambda(6\lambda+120)+3(\lambda^2-\lambda_n^2)^2
 &&\quad\text{(by \eqref{eq:middle-mass-24})}\\
 &<2\lambda(6\lambda+120)+3\cdot30^2
 &&\quad\text{(by \eqref{eq:extreme-square-difference})}\\
 &=12\lambda^2+240\lambda+2700.
\end{alignedat}
\]
Dividing by the positive number $3\lambda^2$ gives
\[
\begin{alignedat}{2}
 S(G)&=\frac{3\lambda^2S(G)}{3\lambda^2}
 &&\quad\text{(by $3\lambda^2>0$)}\\
 &<\frac{12\lambda^2+240\lambda+2700}{3\lambda^2}
 &&\quad\text{(by the preceding chain)}\\
 &=\frac{12\lambda^2}{3\lambda^2}
       +\frac{240\lambda}{3\lambda^2}
       +\frac{2700}{3\lambda^2}\\
 &=4+\frac{80}{\lambda}+\frac{900}{\lambda^2}\\
 &\le4+\frac{80}{1000}+\frac{900}{1000^2}
 &&\quad\text{(by $\lambda\ge1000$)}\\
 &<\frac{41}{10}.
\end{alignedat}
\]
This finishes the proof.
\end{proof}

\begin{lemma}[A second-order perturbation estimate]
\label{lem:perturb}
Let $M_0,E$ be real symmetric matrices of the same order.
Suppose that $v$ is a unit column,
$M_0v=Lv$, and $M_0$ has nonpositive quadratic form on $v^\perp$.
If $\|E\|_{\op}\le b$, $b\ge0$, and $L>2b$, then
\begin{equation}
\label{eq:perturb}
 \lambda_1(M_0+E)\le L+v^{\mathsf T}Ev+
                    \frac{\|Ev\|^2}{L-2b}.
\end{equation}
\end{lemma}

\begin{proof}
In dimension one the only eigenvalue is $L+v^{\mathsf T}Ev$,
so assume the dimension is at least two. Let $V$ be the matrix
whose columns form an orthonormal basis of $v^\perp$. Then
$Q=[\,v\ \ V\,]$ is orthogonal, with $v$ as its first column, and
\begin{equation}
\label{eq:perturb-blocks}
 M_0+E=Q
 \begin{pmatrix}L+v^{\mathsf T}Ev&z^{\mathsf T}\\z&B\end{pmatrix}
 Q^{\mathsf T},
\end{equation}
where $z=V^{\mathsf T}Ev$ and $B=V^{\mathsf T}(M_0+E)V$.
Choose a unit eigenvector $w$ of $B$ for $\lambda_1(B)$.
Since $Vw\in v^\perp$ and $\|Vw\|=\|w\|=1$, we have
\[
\begin{alignedat}{2}
 \lambda_1(B)&=w^{\mathsf T}Bw\\
 &=(Vw)^{\mathsf T}M_0(Vw)+(Vw)^{\mathsf T}E(Vw)\\
 &\le (Vw)^{\mathsf T}E(Vw)
 &&\quad\text{(by nonpositivity of $M_0$ on $v^\perp$)}\\
 &\le \|E\|_{\op}\|Vw\|^2\le b
 &&\quad\text{(by $\|E\|_{\op}\le b$ and $\|Vw\|=1$)}.
\end{alignedat}
\]
The Rayleigh principle gives
\[
\begin{alignedat}{2}
 \lambda_1(M_0+E)&\ge v^{\mathsf T}(M_0+E)v
 &&\quad\text{(by Rayleigh principle)}\\
 &=L+v^{\mathsf T}Ev
 &&\quad\text{(by $M_0v=Lv$, $\|v\|=1$)}\\
 &\ge L-b
 &&\quad\text{(by $\|E\|_{\op}\le b$)}\\
 &>b
 &&\quad\text{(by $L>2b$)}.
\end{alignedat}
\]
Write an eigenvector for $\lambda_1(M_0+E)$ as $Q\binom{\alpha}{y}$,
where $\alpha$ is its coordinate along $v$ and $y$ is the column
of its coordinates in the basis given by $V$. By
\eqref{eq:perturb-blocks} and the eigenvector equation,
\[
 \begin{pmatrix}L+v^{\mathsf T}Ev&z^{\mathsf T}\\ z&B\end{pmatrix}
 \binom{\alpha}{y}=\lambda_1(M_0+E)\binom{\alpha}{y}.
\]
Hence
\begin{equation}
\label{eq:perturb-equations}
 (\lambda_1(M_0+E)-L-v^{\mathsf T}Ev)\alpha=z^{\mathsf T}y,
 \qquad (\lambda_1(M_0+E)I-B)y=\alpha z.
\end{equation}
If $\alpha=0$, the second equation gives $By=\lambda_1(M_0+E)y$
with $y\ne0$, contradicting $\lambda_1(B)\le b<\lambda_1(M_0+E)$.
Thus $\alpha\ne0$.
Since $\lambda_1(M_0+E)>b\ge\lambda_1(B)$, the matrix
$\lambda_1(M_0+E)I-B$ is positive definite. By the spectral theorem,
\[
\begin{alignedat}{2}
 \bigl\|(\lambda_1(M_0+E)I-B)^{-1}\bigr\|_{\op}
 &=\frac{1}{\lambda_{\min}(\lambda_1(M_0+E)I-B)}\\
 &=\frac{1}{\lambda_1(M_0+E)-\lambda_1(B)}\\
 &\le\frac{1}{\lambda_1(M_0+E)-b}
 &&\quad\text{(by $\lambda_1(B)\le b$)}.
\end{alignedat}
\]
Solving the second equation in \eqref{eq:perturb-equations}
for the coordinates in $v^\perp$ and substituting into the first,
we obtain
\[
\begin{alignedat}{2}
 \lambda_1(M_0+E)&=L+v^{\mathsf T}Ev+\frac{z^{\mathsf T}y}{\alpha}
 &&\quad\text{(by the first equation)}\\
 &=L+v^{\mathsf T}Ev+
   \frac{z^{\mathsf T}(\lambda_1(M_0+E)I-B)^{-1}(\alpha z)}{\alpha}
 &&\quad\text{(by the second equation)}\\
 &=L+v^{\mathsf T}Ev+z^{\mathsf T}(\lambda_1(M_0+E)I-B)^{-1}z
 &&\quad\text{(by $\alpha\ne0$)}\\
 &\le L+v^{\mathsf T}Ev+\frac{\|z\|^2}{\lambda_1(M_0+E)-b}
 &&\quad\text{(by the inverse norm bound)}\\
 &\le L+v^{\mathsf T}Ev+\frac{\|Ev\|^2}{\lambda_1(M_0+E)-b}
 &&\quad\text{(by orthogonal projection)}\\
 &\le L+v^{\mathsf T}Ev+\frac{\|Ev\|^2}{L-2b}
 &&\quad\text{(by $\lambda_1(M_0+E)\ge L-b$)}.
\end{alignedat}
\]
This finishes the proof.
\end{proof}

\begin{lemma}[A rank-two matrix with a bounded perturbation]
\label{lem:unified}
Let $a,f$ be orthogonal unit column vectors in $\R^r$, let $L\ge10000$,
and let $E$ be a real symmetric $r\times r$ matrix. Suppose that
\begin{equation}
\label{eq:core-hyp}
 A=L(af^{\mathsf T}+fa^{\mathsf T})+E,\qquad
 \|E\|_{\op}\le7,\qquad Ef=0,\qquad a^{\mathsf T}Ea=0,
\end{equation}
and
\begin{equation}
\label{eq:core-energy}
 \|E\|_{\Frob}^2-2\|Ea\|^2\ge5.
\end{equation}
Then $S(A)>9/2$.
\end{lemma}

\begin{proof}
Write $M_0=L(af^{\mathsf T}+fa^{\mathsf T})$. The two unit columns
$(a+f)/\sqrt2$ and $(a-f)/\sqrt2$ satisfy
\[
 M_0\frac{a\pm f}{\sqrt2}=\pm L\frac{a\pm f}{\sqrt2},
\]
and $M_0$ is zero on $\operatorname{span}\{a,f\}^{\perp}$.
For these columns,
\[
\begin{alignedat}{2}
 \left(\frac{a\pm f}{\sqrt2}\right)^{\mathsf T}
 E\left(\frac{a\pm f}{\sqrt2}\right)
 &=\frac{a^{\mathsf T}Ea\pm a^{\mathsf T}Ef
               \pm f^{\mathsf T}Ea+f^{\mathsf T}Ef}{2}\\
 &=0
 &&\quad\text{(by \eqref{eq:core-hyp} and symmetry)}.
\end{alignedat}
\]
Also,
\[
\begin{alignedat}{2}
 \left\|E\frac{a\pm f}{\sqrt2}\right\|^2
 &=\frac{\|Ea\pm Ef\|^2}{2}\\
 &=\frac{\|Ea\|^2}{2}
 &&\quad\text{(by $Ef=0$)}\\
 &\le\frac{\|E\|_{\op}^2\|a\|^2}{2}
 &&\quad\text{(by the operator norm inequality)}\\
 &\le\frac{49}{2}
 &&\quad\text{(by $\|E\|_{\op}\le7$, $\|a\|=1$)}.
\end{alignedat}
\]
Apply Lemma~\ref{lem:perturb} to $M_0+E$ with $(a+f)/\sqrt2$,
and to $-M_0-E$ with $(a-f)/\sqrt2$. The resulting bounds are
\begin{equation}
\label{eq:extreme-core}
 \lambda_1(A)\le L+\frac{\|Ea\|^2}{2(L-14)},\qquad
 -\lambda_r(A)\le L+\frac{\|Ea\|^2}{2(L-14)}.
\end{equation}
The Rayleigh quotients of these same two columns for $A$ are
$L$ and $-L$, respectively. Hence $\lambda_1(A)\ge L>0$ and
$-\lambda_r(A)\ge L>0$, so the quantities bounded in
\eqref{eq:extreme-core} can be squared.

And we have
\[
\begin{alignedat}{2}
 \tr(M_0E)&=L(f^{\mathsf T}Ea+a^{\mathsf T}Ef)
 &&\quad\text{(by the definition of $M_0$)}\\
 &=0
 &&\quad\text{(by $Ef=0$, $E$ symmetric)},\\
 \tr(M_0^2)&=L^2+(-L)^2
 &&\quad\text{(by the nonzero eigenvalues of $M_0$)}\\
 &=2L^2.
\end{alignedat}
\]
Therefore
\begin{equation}
\label{eq:core-trace}
\begin{alignedat}{2}
 \tr(A^2)&=\tr((M_0+E)^2)
 &&\quad\text{(by $A=M_0+E$)}\\
 &=\tr(M_0^2)+2\tr(M_0E)+\tr(E^2)\\
 &=2L^2+\|E\|_{\Frob}^2
 &&\quad\text{(by the trace values above)}.
\end{alignedat}
\end{equation}
Subtracting the two extreme squares now gives
\begin{equation}
\label{eq:core-mass-bound}
\begin{alignedat}{2}
 S(A)&=\tr(A^2)-\lambda_1(A)^2-\lambda_r(A)^2
 &&\quad\text{(by \eqref{eq:middle-mass})}\\
 &=2L^2+\|E\|_{\Frob}^2-\lambda_1(A)^2-\lambda_r(A)^2
 &&\quad\text{(by \eqref{eq:core-trace})}\\
 &\ge2L^2+\|E\|_{\Frob}^2
        -2\left(L+\frac{\|Ea\|^2}{2(L-14)}\right)^2
 &&\quad\text{(by \eqref{eq:extreme-core})}\\
 &=\|E\|_{\Frob}^2-2\|Ea\|^2
        -\frac{28\|Ea\|^2}{L-14}
        -\frac{\|Ea\|^4}{2(L-14)^2}.
\end{alignedat}
\end{equation}
Since $L-14\ge9986$ and $\|Ea\|^2\le49<50$, the error satisfies
\begin{equation}
\label{eq:core-error}
\begin{alignedat}{2}
 &\frac{28\|Ea\|^2}{L-14}+\frac{\|Ea\|^4}{2(L-14)^2}\\
 &\le\frac{1400}{L-14}+\frac{1250}{(L-14)^2}
 &&\quad\text{(by $\|Ea\|^2<50$)}\\
 &\le\frac{1400}{9986}+\frac{1250}{9986^2}
 &&\quad\text{(by $L-14\ge9986$)}\\
 &<\frac12.
\end{alignedat}
\end{equation}
Combining the last two estimates yields
\[
\begin{alignedat}{2}
 S(A)&\ge\|E\|_{\Frob}^2-2\|Ea\|^2
       -\frac{28\|Ea\|^2}{L-14}-\frac{\|Ea\|^4}{2(L-14)^2}
 &&\quad\text{(by \eqref{eq:core-mass-bound})}\\
 &>\|E\|_{\Frob}^2-2\|Ea\|^2-\frac12
 &&\quad\text{(by \eqref{eq:core-error})}\\
 &\ge5-\frac12
 &&\quad\text{(by \eqref{eq:core-energy})}\\
 &=\frac92.
\end{alignedat}
\]
This finishes the proof.
\end{proof}

\section{Two induced configurations}
\label{sec:configurations}

We next estimate the middle eigenvalues of two induced
configurations. They will exclude three edges in the neighbourhood
of a Perron maximum and an edge between its nonneighbours.

In this section fix
\begin{equation}
\label{eq:reservoir-size}
 N=10^8.
\end{equation}
The role of this size is to ensure that the rank-two term in
Lemma~\ref{lem:unified} has $L\ge\sqrt N=10000$.

\begin{lemma}[A vertex joined to a small graph and an independent set]
\label{lem:cone}
Let $J$ have the vertex partition $\{u\}\sqcup Z\sqcup X$, where
$|Z|\le6$, $e(J[Z])\ge3$, and $X$ is an independent set of size
$N$. Suppose that $u$ is adjacent to every vertex of $Z\cup X$
and that there are no edges between $Z$ and $X$. Then $S(J)>9/2$.
\end{lemma}

\begin{proof}
Order the coordinates as $u,Z,X$. The adjacency matrix is
\begin{equation}
\label{eq:cone-adjacency}
 A_J=\begin{pmatrix}
 0&\one^{\mathsf T}&\one^{\mathsf T}\\
 \one&A_Z&0\\
 \one&0&0
 \end{pmatrix}.
\end{equation}
In the same block order, take
\[
 a=\begin{pmatrix}1\\0\\0\end{pmatrix},
 \qquad f=\frac1{\sqrt N}\begin{pmatrix}0\\0\\\one\end{pmatrix},
 \qquad L=\sqrt N=10000.
\]
These are orthogonal unit columns, and
\begin{equation}
\label{eq:cone-decomposition}
 A_J=L(af^{\mathsf T}+fa^{\mathsf T})+E,
 \qquad
 E=\begin{pmatrix}
 0&\one^{\mathsf T}&0\\
 \one&A_Z&0\\
 0&0&0
 \end{pmatrix}.
\end{equation}
Multiplying by the two block vectors gives
\[
 Ea=E\begin{pmatrix}1\\0\\0\end{pmatrix}
    =\begin{pmatrix}0\\\one\\0\end{pmatrix},
 \qquad
 Ef=\frac1{\sqrt N}E\begin{pmatrix}0\\0\\\one\end{pmatrix}
    =\begin{pmatrix}0\\0\\0\end{pmatrix}.
\]
Thus $a^{\mathsf T}Ea=0$ and $\|Ea\|^2=|Z|$.
Each edge represented in $E$ contributes two entries equal to $1$.
Since $e(J[Z])\le\binom62=15$,
\[
\begin{alignedat}{2}
 \|E\|_{\Frob}^2&=2\bigl(|Z|+e(J[Z])\bigr)
 &&\quad\text{(by counting the edges)}\\
 &\le2(6+15)=42<49
 &&\quad\text{(by $|Z|\le6$, $e(J[Z])\le15$)},\\
 \|E\|_{\Frob}^2-2\|Ea\|^2
 &=2\bigl(|Z|+e(J[Z])\bigr)-2|Z|\\
 &=2e(J[Z])\\
 &\ge6
 &&\quad\text{(by $e(J[Z])\ge3$)}.
\end{alignedat}
\]
Also,
\[
 \|E\|_{\op}\le\|E\|_{\Frob}\le\sqrt{42}<7.
\]
Thus \eqref{eq:cone-decomposition} satisfies every hypothesis of
Lemma~\ref{lem:unified}, including its energy bound, and $S(J)>9/2$.
This finishes the proof.
\end{proof}

\begin{lemma}[An edge between two nonneighbours]
\label{lem:nonneighbor-core}
Let $J$ have the vertex partition $\{u,w,z\}\sqcup Z\sqcup X$,
where $|Z|\le4$, $e(J[Z])=2$, and $X$ is an independent set of
size $N$. Suppose that $u$ is adjacent to every vertex of $Z\cup X$,
that $uw,uz$ are nonedges, that $wz$ is an edge, and that there
are no edges between $Z$ and $X$. Suppose also that one of the
following holds throughout $X$: no vertex of $X$ is adjacent to
$w$ or $z$; every vertex of $X$ is adjacent to $w$ and not $z$;
or every vertex of $X$ is adjacent to $z$ and not $w$.
Then $S(J)>9/2$.
\end{lemma}

\begin{proof}
Let $b,c\in\{0,1\}^{Z}$ record adjacency from $w,z$ to $Z$:
$b_i=1$ precisely when $wi\in E(J)$, and $c_i=1$ precisely when
$zi\in E(J)$. In particular,
\[
 \sum_{i\in Z}b_i=d_Z(w),\qquad
 \sum_{i\in Z}c_i=d_Z(z).
\]
On the ordered coordinates $u,w,z,Z$, the adjacency matrix is
\begin{equation}
\label{eq:B-core}
 B=\begin{pmatrix}
 0&0&0&\one^{\mathsf T}\\
 0&0&1&b^{\mathsf T}\\
 0&1&0&c^{\mathsf T}\\
 \one&b&c&A_Z
 \end{pmatrix}.
\end{equation}
Its edges are those from $u,w,z$ to $Z$, the edge $wz$, and the
two edges in $J[Z]$. Hence
\begin{equation}
\label{eq:B-frobenius}
\begin{alignedat}{2}
 \|B\|_{\Frob}^2
 &=2\bigl(|Z|+d_Z(w)+d_Z(z)+1+2\bigr)
 &&\quad\text{(by counting the edges)}\\
 &\le2\bigl(|Z|+|Z|+|Z|+3\bigr)
 &&\quad\text{(by $d_Z(w),d_Z(z)\le|Z|$)}\\
 &\le2(4+4+4+3)=30
 &&\quad\text{(by $|Z|\le4$)}.
\end{alignedat}
\end{equation}

According to the three possible neighbourhoods of $X$, let $r$
be the column on $\{u,w,z\}\cup Z$ given respectively by
\begin{equation}
\label{eq:attachment-columns}
 r=\begin{pmatrix}1\\0\\0\\0\end{pmatrix},
 \qquad r=\begin{pmatrix}1\\1\\0\\0\end{pmatrix},
 \qquad r=\begin{pmatrix}1\\0\\1\\0\end{pmatrix}.
\end{equation}
Every vertex of $X$ has this same adjacency column. Ordering the
full vertex set as $u,w,z,Z,X$ gives
\begin{equation}
\label{eq:nonneighbor-adjacency}
 A_J=\begin{pmatrix}B&r\one^{\mathsf T}\\
                     \one r^{\mathsf T}&0\end{pmatrix}.
\end{equation}
In the same two blocks as \eqref{eq:nonneighbor-adjacency}, take
\begin{equation}
\label{eq:nonneighbor-decomposition}
 a=\frac1{\|r\|}\binom r0,\qquad
 f=\frac1{\sqrt N}\binom0{\one},\qquad
 L=\sqrt N\|r\|,\qquad E=\begin{pmatrix}B&0\\0&0\end{pmatrix}.
\end{equation}
The block product in \eqref{eq:nonneighbor-adjacency} gives
$A_J=L(af^{\mathsf T}+fa^{\mathsf T})+E$. The supports of $a,f$
are disjoint, so they are orthogonal unit columns. Also $Ef=0$,
and $L\ge10000$, since $\|r\|$ is $1$ or $\sqrt2$.
The support of each column in \eqref{eq:attachment-columns} spans
no edge, because $uw,uz$ are nonedges. Thus
\[
 a^{\mathsf T}Ea=\frac{r^{\mathsf T}Br}{\|r\|^2}=0.
\]
Moreover,
\[
\begin{alignedat}{2}
 \|E\|_{\op}&\le\|E\|_{\Frob}
 &&\quad\text{(by Cauchy--Schwarz)}\\
 &=\|B\|_{\Frob}\\
 &\le\sqrt{30}<7
 &&\quad\text{(by \eqref{eq:B-frobenius})}.
\end{alignedat}
\]
It remains to check $\|E\|_{\Frob}^2-2\|Ea\|^2\ge5$ in each
of the three cases.
In the coordinate order $u,w,z,Z,X$, the block multiplication is
\[
 Ea=\begin{pmatrix}B&0\\0&0\end{pmatrix}
       \binom{r/\|r\|}{0}
    =\binom{Br/\|r\|}{0}.
\]

If the vertices of $X$ are adjacent to neither $w$ nor $z$,
multiplication in \eqref{eq:B-core} gives
\[
 Ea=\begin{pmatrix}0\\0\\0\\\one\\0\end{pmatrix},
 \qquad \|Ea\|^2=|Z|.
\]
Therefore
\[
\begin{alignedat}{2}
 &\|E\|_{\Frob}^2-2\|Ea\|^2\\
 &=2\bigl(|Z|+d_Z(w)+d_Z(z)+3\bigr)-2|Z|\\
 &=6+2d_Z(w)+2d_Z(z)\\
 &\ge6
 &&\quad\text{(by nonnegativity of degrees)}.
\end{alignedat}
\]

If the vertices of $X$ are adjacent to $w$ and not $z$, then
\[
 Ea=\frac1{\sqrt2}\begin{pmatrix}0\\0\\1\\\one+b\\0\end{pmatrix}.
\]
Consequently,
\[
\begin{alignedat}{2}
 \|Ea\|^2&=\frac{1+\sum_{i\in Z}(1+b_i)^2}{2}\\
 &=\frac{1+\sum_{i\in Z}(1+3b_i)}{2}
 &&\quad\text{(by $b_i^2=b_i$)}\\
 &=\frac{1+|Z|+3\sum_{i\in Z}b_i}{2}\\
 &=\frac{|Z|+3d_Z(w)+1}{2}
 &&\quad\text{(by $\sum_{i\in Z}b_i=d_Z(w)$)}.
\end{alignedat}
\]
It follows that
\[
\begin{alignedat}{2}
 &\|E\|_{\Frob}^2-2\|Ea\|^2\\
 &=2\bigl(|Z|+d_Z(w)+d_Z(z)+3\bigr)
       -\bigl(|Z|+3d_Z(w)+1\bigr)\\
 &=5+|Z|-d_Z(w)+2d_Z(z)\\
 &\ge5
 &&\quad\text{(by $d_Z(w)\le|Z|$, $d_Z(z)\ge0$)}.
\end{alignedat}
\]

The case in which the vertices of $X$ are adjacent to $z$ and
not $w$ follows by interchanging $w$ and $z$.
Every case satisfies \eqref{eq:core-energy}, so
Lemma~\ref{lem:unified} gives $S(J)>9/2$.
This finishes the proof.
\end{proof}

\section{The neighbourhood of a Perron maximum}
\label{sec:neighbourhood}

Let $G$ be a simple graph on $n$ vertices with $m$ edges satisfying
$\lambda=\lambda(G)>0$. By the Perron--Frobenius theorem, choose
a nonnegative Perron vector $x=(x_v)_{v\in V(G)}$ and a vertex
$u$ corresponding to a largest entry of $x$. Normalize $x$ so that
\begin{equation}
\label{eq:perron}
 Ax=\lambda x,\qquad x_u=1,\qquad0\le x_v\le1
 \quad(v\in V(G)).
\end{equation}
Partition the vertices into
\begin{equation}
\label{eq:partition}
 V(G)=\{u\}\sqcup U\sqcup W,\qquad
 U=N(u),\qquad W=V(G)\setminus(\{u\}\cup N(u)).
\end{equation}
The equation at $u$ gives
\begin{equation}
\label{eq:root-eigen}
\begin{alignedat}{2}
 \lambda&=\sum_{v\in U}x_v
 &&\quad\text{(by $x_u=1$, $U=N(u)$)}\\
 &\le\sum_{v\in U}1
 &&\quad\text{(by $x_v\le1$)}\\
 &=|U|.
\end{alignedat}
\end{equation}
Write
\begin{equation}
\label{eq:weighted-sums}
 F_x=\sum_{ab\in E(G[U])}(x_a+x_b),\qquad
 M_x=\sum_{y\in W}d_U(y)(1-x_y).
\end{equation}
Each unordered edge occurs once in $F_x$. Thus
$0\le F_x\le2e(G[U])$ and $M_x\ge0$.

\begin{lemma}[The two-step equation and the edge count]
\label{lem:cut}
Let $G$ be a simple graph on $n$ vertices with $m$ edges and
$\lambda=\lambda(G)>0$. With the Perron vector and vertex partition
\eqref{eq:perron}--\eqref{eq:partition}, one has $e(G[U])\le t=t(G)$ and
\begin{equation}
\label{eq:cut-identity}
 \lambda^2-m=F_x-e(G[U])-e(G[W])-M_x.
\end{equation}
If $\lambda^2>m$, then $e(G[W])<e(G[U])$.
\end{lemma}

\begin{proof}
An edge $ab\in E(G[U])$ gives the triangle $\{u,a,b\}$, and
deleting $u$ recovers the edge. This is an injection into the
triangles of $G$, so
\begin{equation}
\label{eq:neighbourhood-triangles}
 e(G[U])\le t.
\end{equation}
The edges of $G$ consist of the $|U|$ edges at $u$, the edges in
$U$, the edges between $U$ and $W$, and the edges in $W$. Therefore
\begin{equation}
\label{eq:edge-partition}
 m=|U|+e(G[U])+\sum_{y\in W}d_U(y)+e(G[W]).
\end{equation}
Sum the eigenvector equations over $U$ and then group their
neighbours according to \eqref{eq:partition}. This gives
\begin{equation}
\label{eq:root-two-step}
\begin{alignedat}{2}
 \lambda^2&=\lambda\sum_{v\in U}x_v
 &&\quad\text{(by \eqref{eq:root-eigen})}\\
 &=\sum_{v\in U}\sum_{z\in N(v)}x_z
 &&\quad\text{(by the eigenvector equations)}\\
 &=|U|+\sum_{ab\in E(G[U])}(x_a+x_b)
                  +\sum_{y\in W}d_U(y)x_y\\
 &=|U|+F_x+\sum_{y\in W}d_U(y)x_y
 &&\quad\text{(by \eqref{eq:weighted-sums})}.
\end{alignedat}
\end{equation}

Subtracting \eqref{eq:edge-partition} from \eqref{eq:root-two-step}
now yields
\[
\begin{alignedat}{2}
 \lambda^2-m
 &=|U|+F_x+\sum_{y\in W}d_U(y)x_y-m
 &&\quad\text{(by \eqref{eq:root-two-step})}\\
 &=F_x+\sum_{y\in W}d_U(y)x_y-e(G[U])
           -\sum_{y\in W}d_U(y)-e(G[W])
 &&\quad\text{(by \eqref{eq:edge-partition})}\\
 &=F_x-e(G[U])-e(G[W])
           -\sum_{y\in W}d_U(y)(1-x_y)\\
 &=F_x-e(G[U])-e(G[W])-M_x
 &&\quad\text{(by \eqref{eq:weighted-sums})}.
\end{alignedat}
\]
This proves \eqref{eq:cut-identity}. If $\lambda^2>m$, it also gives
\begin{equation}
\label{eq:W-edges-less-U}
\begin{alignedat}{2}
 e(G[W])&<F_x-e(G[U])-M_x
 &&\quad\text{(by \eqref{eq:cut-identity})}\\
 &\le2e(G[U])-e(G[U])-M_x
 &&\quad\text{(by $F_x\le2e(G[U])$)}\\
 &=e(G[U])-M_x\\
 &\le e(G[U])
 &&\quad\text{(by $M_x\ge0$)}.
\end{alignedat}
\end{equation}
This finishes the proof.
\end{proof}

For the rest of the proof, fix the absolute constants
\begin{equation}
\label{eq:constants}
 R=\binom{N+5}{6},\qquad C=R+10,\qquad N=10^8.
\end{equation}
Clearly, $R>3N$ and $C>1000$.

The classical Ramsey bound of Erd\H{o}s and
Szekeres~\cite{ErdosSzekeres1935} says that a graph with at least
$\binom{s+q-2}{s-1}$ vertices contains $K_s$ or an independent
set of size $q$. For $s=7$ and $q=N$, this means that every
$K_7$-free graph on at least $R$ vertices contains an independent
set of size $N$.

\begin{lemma}[The neighbourhood structure]
\label{lem:structure}
Let $G$ be a simple graph on $n$ vertices with $m$ edges satisfying
\[
 \lambda=\lambda(G)>\sqrt m,\qquad S(G)<41/10,\qquad t=t(G)<\lambda-C,
\]
and suppose that $G$ has no common triangle vertex. Choose $u,x$
as in \eqref{eq:perron}, and use \eqref{eq:partition}. Then $W$
is independent, and $G[U]$ consists of two vertex-disjoint edges
and isolated vertices.
\end{lemma}

\begin{proof}
Lemma~\ref{lem:interlacing} shows that $G$ contains no $K_7$.
We first prove that $G[U]$ has exactly two edges.

Suppose that $e(G[U])\ge3$. Choose three distinct edges in $G[U]$
and let $Z$ be the set of their endpoints. Then
\[
 |Z|\le6,\qquad 3\le e(G[Z])\le\binom62=15.
\]
Let $Y$ consist of the vertices in $U\setminus Z$ with no neighbour
in $Z$, and let $B=U\setminus(Z\cup Y)$. Thus
\[
 U=Z\sqcup Y\sqcup B.
\]
For each $b\in B$, choose an edge from $b$ to $Z$. These edges
are distinct because their endpoints outside $Z$ are distinct,
and none is an edge of $G[Z]$. Hence
$|B|\le e(G[U])-e(G[Z])$. Consequently,
\[
\begin{alignedat}{2}
 |Y|&=|U|-|Z|-|B|
 &&\quad\text{(by the vertex partition)}\\
 &\ge|U|-|Z|-e(G[U])+e(G[Z])
 &&\quad\text{(by the bound on $|B|$)}\\
 &\ge|U|-|Z|-t+e(G[Z])
 &&\quad\text{(by \eqref{eq:neighbourhood-triangles})}\\
 &\ge\lambda-t-|Z|+e(G[Z])
 &&\quad\text{(by \eqref{eq:root-eigen})}\\
 &>C-|Z|+e(G[Z])
 &&\quad\text{(by $t<\lambda-C$)}\\
 &\ge C-3
 &&\quad\text{(by $|Z|\le6$, $e(G[Z])\ge3$)}\\
 &=R+7>R
 &&\quad\text{(by \eqref{eq:constants})}.
\end{alignedat}
\]
The Ramsey bound gives an independent set $X\subseteq Y$ of size
$N$. The induced graph $J=G[\{u\}\cup Z\cup X]$ satisfies
Lemma~\ref{lem:cone}: $u$ is adjacent to all of $Z\cup X$, and
there are no edges between $Z$ and $X$. Lemmas~\ref{lem:cone}
and~\ref{lem:interlacing} would give
\[
 \frac92<S(J)\le S(G)<\frac{41}{10},
\]
a contradiction. Therefore $e(G[U])\le2$.

If $e(G[U])=0$, then \eqref{eq:W-edges-less-U} gives
$e(G[W])<0$, which is impossible. If $e(G[U])=1$, the same bound
makes $W$ independent. Let $ab$ be the unique edge in $G[U]$.
Every triangle containing $u$ is $\{u,a,b\}$. A triangle avoiding
$u$ cannot contain two vertices of $W$, since $W$ is independent,
and cannot have all three vertices in $U$, since $G[U]$ has only
one edge. Its two vertices in $U$ must therefore be $a,b$.
Thus every triangle contains $a$, contrary to the hypothesis.
We have proved
\begin{equation}
\label{eq:two-neighbourhood-edges}
 e(G[U])=2.
\end{equation}

We next prove that $W$ is independent. Suppose instead that
$wz\in E(G[W])$. Now let $Z$ be the set of endpoints of the two
edges in $G[U]$. Then $|Z|\le4$ and $e(G[Z])=2$. Put
\[
 T=U\cap N(w)\cap N(z),\qquad Y=U\setminus(Z\cup T).
\]
Each $v\in T$ gives a triangle $\{w,z,v\}$, and these triangles
are distinct. None contains $u$, so they are different from the
two triangles at $u$ supplied by \eqref{eq:two-neighbourhood-edges}.
Therefore
\begin{equation}
\label{eq:T-count}
 |T|+2\le t.
\end{equation}
Every vertex of $Y$ is isolated in $G[U]$, since $Z$ contains all
endpoints of its two edges, and no vertex of $Y$ is adjacent to
both $w$ and $z$. Moreover,
\[
\begin{alignedat}{2}
 |Y|&=|U|-|Z\cup T|
 &&\quad\text{(by $Z,T\subseteq U$)}\\
 &\ge|U|-|Z|-|T|
 &&\quad\text{(by the union bound)}\\
 &\ge|U|-4-|T|
 &&\quad\text{(by $|Z|\le4$)}\\
 &\ge|U|-4-(t-2)
 &&\quad\text{(by \eqref{eq:T-count})}\\
 &=|U|-t-2\\
 &\ge\lambda-t-2
 &&\quad\text{(by \eqref{eq:root-eigen})}\\
 &>C-2=R+8
 &&\quad\text{(by $t<\lambda-C$, \eqref{eq:constants})}\\
 &>3N
 &&\quad\text{(by $R>3N$)}.
\end{alignedat}
\]
Partition $Y$ by its adjacency to $w,z$ as
\[
 Y_0=Y\setminus(N(w)\cup N(z)),\qquad
 Y_w=Y\cap(N(w)\setminus N(z)),\qquad
 Y_z=Y\cap(N(z)\setminus N(w)).
\]
These sets are disjoint and cover $Y$. Since $|Y|>3N$, one of
them contains at least $N$ vertices. Choose $N$ vertices from that set as
$X$. They are independent and have no neighbours in $Z$, because
they are isolated in $G[U]$. The induced graph
$J=G[\{u,w,z\}\cup Z\cup X]$ now satisfies every hypothesis of
Lemma~\ref{lem:nonneighbor-core}. Once again,
\[
 \frac92<S(J)\le S(G)<\frac{41}{10}
\]
is impossible. Thus
\begin{equation}
\label{eq:W-independent}
 e(G[W])=0.
\end{equation}

Finally, the two edges in $G[U]$ cannot share an endpoint $a$.
If they did, every triangle containing $u$ would contain $a$.
A triangle avoiding $u$ would have exactly two vertices in $U$,
since there are no edges in $W$ and a graph with only two edges has no
triangle. Those two vertices span one of the edges in $G[U]$,
which again contains $a$. Every triangle would contain $a$,
contrary to the hypothesis. Therefore the two edges are
vertex-disjoint. This finishes the proof.
\end{proof}

\section{Proof of the main theorem}
\label{sec:main-proof}

For a graph with $\lambda>0$, choose a normalized nonnegative
Perron vector and a vertex $u$ as in \eqref{eq:perron}. Continue
to use $U,W,F_x,M_x$ from \eqref{eq:partition} and
\eqref{eq:weighted-sums}.

\begin{lemma}[The estimate for two disjoint neighbourhood edges]
\label{lem:matching}
Let $G$ be a simple graph on $n$ vertices with $m$ edges satisfying
$\lambda^2>m$, where $\lambda=\lambda(G)$. Suppose that $W$ is
independent and that
\begin{equation}
\label{eq:matching}
 E(G[U])=\{ab,cd\}
\end{equation}
for four distinct vertices $a,b,c,d$. Then $\lambda<t+2$, where
$t=t(G)$.
\end{lemma}

\begin{proof}
Let $D=\{a,b,c,d\}$. The definition of $F_x$ gives
$F_x=x_a+x_b+x_c+x_d\le4$. By Lemma~\ref{lem:cut},
\begin{equation}
\label{eq:F-minus-M}
\begin{alignedat}{2}
 F_x-M_x&=\lambda^2-m+e(G[U])+e(G[W])
 &&\quad\text{(by \eqref{eq:cut-identity})}\\
 &=\lambda^2-m+2
 &&\quad\text{(by \eqref{eq:matching}, $W$ independent)}\\
 &>2
 &&\quad\text{(by $\lambda^2>m$)}.
\end{alignedat}
\end{equation}

For $y\in W$, let $q_y$ count the two matching edges whose
endpoints are both adjacent to $y$:
\begin{equation}
\label{eq:q-definition}
 q_y=\one_{\{a,b\}\subseteq N(y)}+
     \one_{\{c,d\}\subseteq N(y)}.
\end{equation}
Here $\one_P$ is $1$ if the assertion $P$ holds and $0$ otherwise.
The only triangles containing $u$ are $\{u,a,b\}$ and
$\{u,c,d\}$. A triangle avoiding $u$ has exactly one vertex in $W$,
since $W$ is independent and $G[U]$ has only two edges. Its two
vertices in $U$ form either $ab$ or $cd$, and its vertex in $W$
is a common neighbour of the endpoints of that edge. Thus
\begin{equation}
\label{eq:triangle-count}
\begin{aligned}
 t&=2+|W\cap N(a)\cap N(b)|+|W\cap N(c)\cap N(d)|\\
  &=2+\sum_{y\in W}
     \bigl(\one_{\{a,b\}\subseteq N(y)}+
           \one_{\{c,d\}\subseteq N(y)}\bigr)\\
  &=2+\sum_{y\in W}q_y.
\end{aligned}
\end{equation}

The equations at the four endpoints are
\begin{equation}
\label{eq:four-equations}
\begin{aligned}
 \lambda x_a&=1+x_b+\sum_{y\in W\cap N(a)}x_y,&
 \lambda x_b&=1+x_a+\sum_{y\in W\cap N(b)}x_y,\\
 \lambda x_c&=1+x_d+\sum_{y\in W\cap N(c)}x_y,&
 \lambda x_d&=1+x_c+\sum_{y\in W\cap N(d)}x_y.
\end{aligned}
\end{equation}
Adding these equations and collecting terms for each $y$ gives
\begin{equation}
\label{eq:four-sum}
\begin{alignedat}{2}
 \lambda F_x&=\sum_{v\in D}\lambda x_v
 &&\quad\text{(by the definition of $F_x$)}\\
 &=4+F_x+\sum_{v\in D}\sum_{y\in W\cap N(v)}x_y
 &&\quad\text{(by \eqref{eq:four-equations})}\\
 &=4+F_x+\sum_{y\in W}\sum_{v\in D\cap N(y)}x_y\\
 &=4+F_x+\sum_{y\in W}|N(y)\cap D|x_y.
\end{alignedat}
\end{equation}

We next bound the last sum. First, for $v\in D$,
\begin{equation}
\label{eq:endpoint-lower}
\begin{alignedat}{2}
 \lambda x_v&=\sum_{z\in N(v)}x_z
 &&\quad\text{(by the eigenvector equation)}\\
 &\ge x_u
 &&\quad\text{(by $u\in N(v)$, $x_z\ge0$)}\\
 &=1
 &&\quad\text{(by normalization)}.
\end{alignedat}
\end{equation}
For $y\in W$, independence of $W$ and the nonedge $uy$ imply
$N(y)\subseteq U$. Hence
\begin{equation}
\label{eq:y-upper}
\begin{alignedat}{2}
 \lambda x_y&=\sum_{z\in U\cap N(y)}x_z
 &&\quad\text{(by the eigenvector equation)}\\
 &\le\sum_{z\in U\cap N(y)}1
 &&\quad\text{(by $x_z\le1$)}\\
 &=d_U(y).
\end{alignedat}
\end{equation}
For every $y\in W$ and $v\in D$, these bounds give
\begin{equation}
\label{eq:missing-endpoint}
\begin{alignedat}{2}
 x_y&\le(\lambda x_v)x_y
 &&\quad\text{(by \eqref{eq:endpoint-lower}, $x_y\ge0$)}\\
 &=(\lambda x_y)x_v\\
 &\le d_U(y)x_v
 &&\quad\text{(by \eqref{eq:y-upper}, $x_v\ge0$)}.
\end{alignedat}
\end{equation}
The total weight of the vertices of $U$ not adjacent to $y$ is
\begin{equation}
\label{eq:missing-U}
\begin{alignedat}{2}
 \sum_{z\in U\setminus N(y)}x_z
 &=\sum_{z\in U}x_z-\sum_{z\in U\cap N(y)}x_z\\
 &=\lambda-\lambda x_y
 &&\quad\text{(by \eqref{eq:root-eigen}, \eqref{eq:y-upper})}\\
 &=\lambda(1-x_y).
\end{alignedat}
\end{equation}

Fix $y\in W$ and one matching edge $v_1v_2\in\{ab,cd\}$. There
are four possibilities for its endpoints:
\[
\begin{array}{c@{\qquad}l@{\qquad}l}
 N(y)\cap\{v_1,v_2\}&\text{inequality}&\text{reason}\\[4pt]
 \{v_1,v_2\}&2x_y\le2&\text{(by $x_y\le1$)},\\[4pt]
 \{v_1\}&x_y\le d_U(y)x_{v_2}&\text{(by \eqref{eq:missing-endpoint})},\\[4pt]
 \{v_2\}&x_y\le d_U(y)x_{v_1}&\text{(by \eqref{eq:missing-endpoint})},\\[4pt]
 \varnothing&0\le d_U(y)(x_{v_1}+x_{v_2})&\text{(by nonnegativity)}.
\end{array}
\]
Apply the relevant row to $ab$ and $cd$, then add. Each adjacent
endpoint contributes $x_y$ on the left. On the right, each edge
with both endpoints adjacent to $y$ contributes $2$, and each
missing endpoint contributes its coordinate multiplied by $d_U(y)$.
Hence
\begin{equation}
\label{eq:local-four}
 |N(y)\cap D|x_y\le2q_y+d_U(y)\sum_{v\in D\setminus N(y)}x_v.
\end{equation}
Therefore
\begin{equation}
\label{eq:local-bound}
\begin{alignedat}{2}
 |N(y)\cap D|x_y
 &\le2q_y+d_U(y)\sum_{v\in D\setminus N(y)}x_v
 &&\quad\text{(by \eqref{eq:local-four})}\\
 &\le2q_y+d_U(y)\sum_{z\in U\setminus N(y)}x_z
 &&\quad\text{(by $D\subseteq U$, $x_z\ge0$)}\\
 &=2q_y+\lambda d_U(y)(1-x_y)
 &&\quad\text{(by \eqref{eq:missing-U})}.
\end{alignedat}
\end{equation}
Summing over $y\in W$, we obtain
\begin{equation}
\label{eq:sum-local}
\begin{alignedat}{2}
 &\sum_{y\in W}|N(y)\cap D|x_y\\
 &\le\sum_{y\in W}[2q_y+\lambda d_U(y)(1-x_y)]
 &&\quad\text{(by \eqref{eq:local-bound})}\\
 &=2\sum_{y\in W}q_y+\lambda\sum_{y\in W}d_U(y)(1-x_y)\\
 &=2\sum_{y\in W}q_y+\lambda M_x
 &&\quad\text{(by \eqref{eq:weighted-sums})}.
\end{alignedat}
\end{equation}
Substitution into the endpoint sum now gives
\[
\begin{alignedat}{2}
 \lambda F_x&=4+F_x+\sum_{y\in W}|N(y)\cap D|x_y
 &&\quad\text{(by \eqref{eq:four-sum})}\\
 &\le4+F_x+2\sum_{y\in W}q_y+\lambda M_x
 &&\quad\text{(by \eqref{eq:sum-local})}\\
 &=4+F_x+2(t-2)+\lambda M_x
 &&\quad\text{(by \eqref{eq:triangle-count})}\\
 &=F_x+2t+\lambda M_x.
\end{alignedat}
\]
Hence
\begin{equation}
\label{eq:final-weighted}
\begin{alignedat}{2}
 \lambda(F_x-M_x)&=\lambda F_x-\lambda M_x\\
 &\le F_x+2t+\lambda M_x-\lambda M_x
 &&\quad\text{(by the preceding chain)}\\
 &=F_x+2t\\
 &\le4+2t
 &&\quad\text{(by $F_x\le4$)}.
\end{alignedat}
\end{equation}
Finally,
\[
\begin{alignedat}{2}
 \lambda&=\frac{2\lambda}{2}\\
 &<\frac{\lambda(F_x-M_x)}2
 &&\quad\text{(by \eqref{eq:F-minus-M}, $\lambda>0$)}\\
 &\le\frac{4+2t}{2}
 &&\quad\text{(by \eqref{eq:final-weighted})}\\
 &=t+2.
\end{alignedat}
\]
This finishes the proof.
\end{proof}

\begin{proof}[Proof of Theorem~\ref{thm:main}]
Assume $G$ is a simple graph on $n$ vertices with $m$ edges
satisfying $\lambda=\lambda(G)>\sqrt m$, and suppose that no
vertex of $G$ belongs to every triangle. Write $t=t(G)$.
It suffices to prove $t\ge\lambda-C$, since $\lambda>\sqrt m$.
Suppose to the contrary that $t<\lambda-C$. Since $t\ge0$ and
$C>1000$, we have
\[
 \lambda>t+C\ge C>1000,\qquad t<\lambda-C<\lambda.
\]
Lemma~\ref{lem:smallmass} therefore gives $S(G)<41/10$.
Choose a nonnegative Perron column and $u$ as in \eqref{eq:perron}.
Lemma~\ref{lem:structure} shows that $W$ is independent and that
$G[U]$ consists of two disjoint edges and isolated vertices.
Since $\lambda^2>m$, Lemma~\ref{lem:matching} gives $\lambda<t+2$.
On the other hand,
\[
\begin{alignedat}{2}
 \lambda&>t+C
 &&\quad\text{(by the contrary assumption)}\\
 &>t+2
 &&\quad\text{(by $C>2$)},
\end{alignedat}
\]
a contradiction. This finishes the proof.
\end{proof}

\section*{Statement on the use of generative AI tools}

The authors used GPT-6.1-sol in the preparation of this work, including the development of
proof arguments, the checking of intermediate steps, and the drafting and revision of the manuscript
and its LaTeX source. All mathematical results, arguments, and proofs were independently reviewed
and verified by the authors, who take full responsibility for the accuracy and content of this work.

\end{document}